\documentstyle{amsppt}
\magnification=1200
\NoBlackBoxes
\def\ss{\vskip.10in}
\def\ls{\vskip.20in}
\def\C{\Bbb C}
\def\P{\Bbb P}

\def\Z{\Bbb Z} 
\def\opluss{\operatornamewithlimits{\oplus}}
\def\lll{$\underline{\hskip.15in}$}
\baselineskip=16pt

\centerline{\bf On the density of ratios of Chern numbers}
\centerline{\bf of embedded threefolds}
\ss
\centerline{Mei-Chu Chang\footnote"*"{partially supported by NSF Grant No. 
D.M.S 9304580}}
\centerline{Math Dept.}
\centerline{UCR}
\centerline{Riverside, CA  92521}
\ls
Let $X$ be a 3-fold of general type with Chern class
$c_i (X) = c_i (T_X) \in H^{2i} (X, \Z)$.  Then cup product gives
the Chern numbers $c_1^3 (X), c_1 c_2(X), c_3(X)\in \Z$.  One
natural question to ask is:  {\it Which point in $\P^2(\Cal Q)$ does
correspond to the Chern numbers $[c_1^3 (X), c_1c_2(X), c_3(X)]$ for
some 3-fold $X$?}

By using Fermat cover (desingularization of the branched cover of 
$\C\P^3$ over an arrangement of planes), Hunt [H] was able to show
that there are two triangles inside which {\it every} point corresponds
to a 3-fold.  Extending Hunt's idea, Liu [C3, L] obtained a larger
region in the affine chart $c_1c_2 \neq 0$ (LMCN in the chart attached).
In [C2], the author gave an explicit description of the ``SCI zone",
the {\it limit points} of the Chern ratios $(c_1^3/c_1c_2, c_3/c_1c_2)
(X)$ of all the complete intersection 3-folds $X$.  However the only known bound of  the region of all possible Chern numbers, under the assumption that either
$X$ is minimal or the canonical divisor $K_X$ is ample, is $0 \leq c_1^3/
c_1c_2 \leq 8/3$.  (The right inequality is Yau's)
\ss
In this paper, we study the bounds of the Chern ratios of 3-folds in 
$\P^5$.  Let $X \subset \P^5$ be a 3-fold with hyperplane class $H$, 
canonical class $K$.  What we can show is the following

\proclaim{Theorem} The limit points of $\{(x,y) = (c_1^3/c_1c_2, c_3/c_1c_2)
(X)\; |\; X \subset \P^5$ has $H^i K^j > 0\}$ lie on the line segment
$x+y = 2,\; 1 \leq x \leq 2\}$.
\endproclaim

\remark{Remark 1}  In [C1], we showed that the limit points of Chern ratios
of determinantal 3-folds in $\P^5$ are the line segment $x+y = 2,\; 1 \leq
x \leq 17/12$.
\endremark

\remark{Remark 2}  Note that $x<2$ means the third Segre class of the 
cotangent bundle $\Omega_X^1$ is negative.  As Robert Braun commented
that when the degree $d$ goes to infinity, 3-folds in $\P^5$ have small
positive $\Omega^1$.
\endremark

\remark{Remark 3} The line $x+y = 2$ is part of the boundary of the SCI
zone, the limit points of Chern ratios of all complete intersection 3-folds.
In fact, the first piece of the "lower" curve of the boundary is
$x+y = 2,\; 1 \leq x \leq 4/3$.
\endremark

Our approach to the problem is expressing the Chern numbers in terms of the 
intersection numbers $H^iK^j$ (Property (1)), and reducing the problem to 
whether $K^3$ is greater than a linear combination of $d^2, dH^2K$, and 
$HK^2$.  (cf. Proposition 6).  For 3-fold $X$ contained in hypersurfaces 
of fixed degree $s$, when the degree $d$ goes to infinity, $K^3$ (which is
greater than 
$d^4/s^3)$ dominates, and the ratios $(c_1^3/c_1c_2, c_3/c_1c_2)$ approximate
to the line $x+y =2$ arbitrarily.  On the other hand, when $s$ goes to 
infinity, we use the positivity of the third Segre class of $N_{X/\P_5} (-1)$
(Proposition 1(ii)), and reduce the problem further to an inequality 
between $d H^2 K$ and $HK^2$ (cf. Proposition 5), and use the genus
bound.

\subheading{Acknowledgement}  Most of the work was done while the author
was a member in the Institute for Advanced Study. She would like to thank
the Institute for 
its financial support and the hospitality of its faculty and staff.  She
also thanks Robert Braun, Rob Lazarsfeld and Scott Nollet for helpful
communication, and Bruce Chalmers for technical support.

\vfill\eject

Let $X$ be a 3-fold of degree $d$ in $\P^5$.  We use the following notation:

$H$ = hyperplane class

$K$ = canonical class

$Y$ = general hyperplane section

$C$ = general curve section

$s(X)$ = smallest degree of hypersurfaces containing $X$
$$
\align
(x,y) &:= (c_1^3 (X)/c_1c_2 (X), c_3(X)/c_1c_2(X))\\
&= (K^3/(K^3 - \Delta_1), \ \ (K^3 - \Delta_2)/(K^3 - \Delta_1)), \ \ 
\text{where}\\
\Delta_1 &:= (d-15)H^2K - 6HK^2\\
\Delta_2 &:= (6d - 70)d + (2d -51) H^2 K - 12 HK^2.\\
\endalign
$$

The last equality follows from
$$
\align
c_2 (X) &= (15-d) H^2 + 6HK + K^2\\
c_3 (X) &= (6d - 70)d + (2d-51) H^2 K - 12HK^2 - K^3,
\tag1
\endalign
$$
which is the generalized double point formula gotten from the normal-tangent
sequence.

The other properties we need are the following:  (please see [BOSS] as a
general reference.)

(2) (GHIT, generalized Hodge index theorem) $d HK^2 \leq
(H^2 K)^2$

(3)  (special genus formula) $2g - 2 \leq d^2/s + d (s-4)$

(4) (Castelnuovo)
 $$\aligned &p_g (X) \leq 2 {\binom M4}   + {\binom M3}\\
& p_g (Y)\leq 2 {\binom M3} + {\binom M2},\quad \text{where}\ \ M = 
[(d-1)/2]
\endaligned$$

(5) [BOSS] $- 24 \chi(\Cal O_X) \geq d^4/s^3 + l.t.$ in $\sqrt{d}.$

(6)  $\chi (\Cal O_X) = c_1 c_2/24$, combined with property (1), this is

(6$'$) $K^3 = (d-15) H^2K - 6HK^2 - 24 \chi (\Cal O_X)$

\proclaim{Proposition 1} The following hold

(i) $9H^2K + HK^2 \geq d(d-21)$ 

(ii)  $K^3 \geq 8d(d-13) + 2(d-33) H^2K -14HK^2$
\endproclaim

\demo{Proof}  Let $N$ be the normal bundle of $X$ in $\P^5$.  Then
$N(-1)$ is generated by global sections and has $c_1(N(-1)) = 4 H+K$ and 
$c_2(N(-1))=(d-5)H^2 - HK$.  Now (i) is $(c_1^2 - c_2)N(-1)\cdot
H \geq 0$, and (ii)  is $(c_1^3 - 2c_1 c_2) N (-1) \geq 0$.  \qed
\enddemo

\remark{Notation 1}  l.t. means terms of lower degrees (possibly fractional)
in $d$.
\endremark

\proclaim{Proposition 2} The following hold

(i)  $H^2K < d^2/s + l.t.$  

(ii) $HK^2 < d^3/s^2 + l.t.$  

(iii)  $K^3 > d^4/s^3 + l.t.$
\endproclaim

\demo{Proof}  (i) follows from the adjunction formula and Property (3), the 
genus formula.  (ii) follows from (i) and Property (2), GHIT.  Applying
(i), (ii) and Property (5) to Property (6$'$), we have (iii).  \qed
\enddemo

\proclaim{Proposition 3}  for any integer $d_0$, there are finitely many points
$(x,y)$ corresponding to $(c_1^3/c_1c_2, c_3/c_1c_2) (X)$ such that
$\deg X \leq d_0$.
\endproclaim

\demo{Proof}  for each degree $d$ fixed, (i) and (ii) in the previous proposition imply that $H^2K$ and $HK^2$ have only finitely many possible values.  
The finiteness of $K^3$ follows from Propoer (6$'$).  (The finiteness of 
$\chi (\Cal O_X)$ follows from Property (4) and that $h^1 (\Cal O_X)=0$).
In Property (1), we see that the Chern ratios are functions of $H^iK^j$.\qed
\enddemo

Since we are interested in the limit points of the Chern ratios, in our 
statement we are free to exclude finitely many of them.  By Proposition 1,
we are free to exclude finitely many degrees in our inequalities.  So we
use the following

\remark{Notation 2}  $A<.\ \ B$ means that except for 3-folds of finitely many
possible degrees, $A$ is less than $B$.
\endremark

\proclaim{Proposition 4}   $d > cs^{5/3}$, where $c>0$ is a constant
independent of the 3-fold.
\endproclaim

\demo{Proof}  Let $k = s-1$, then $H^0 (\Cal I_X (k))=0$ and we have
$$
\align
h^0 (\Cal O_{\P_5} (k)) &\leq h^0 (\Cal O_X (k))\\
&\leq \chi (\Cal O_X (k)) + h^1 (\Cal O_X (k)) + h^3 (\Cal O_X (k))
\tag a
\endalign
$$
Riemann-Rock gives, for $d\geq 15$
$$
\chi (\Cal O_X (k)) \leq k^3 d/6 + kH^2K/2 + kHK^2/6 + h^0 (\Cal O_X) - 
h^1 (\Cal O_X) +h^2(\Cal O_X) - h^3(\Cal O_X) . 
\tag b
$$
The $H^3$-cohomologies of the sequence
$$
0 \to \Cal O_X (t-1) \to \Cal O_X (t) \to \Cal O_Y (t) \to 0 
\tag*
$$
give that
$$
h^3 (\Cal O_X (k)) \leq h^3 (\Cal O_X). 
\tag c
$$
Similarly, the $H^1$-cohomologies of the same sequence $(*)$ and its
restriction $0\to \Cal O_Y (t-1) \to \Cal O_Y (t) \to \Cal O_C (t)
\to 0$ give
$$
h^1 (\Cal O_X(k)) \leq h^1 (\Cal O_X) + \sum\limits_{t=0}^k
(k+1 -t) h^1 (\Cal O_C (t)).
\tag d
$$
A rough estimate by the genus bound is
$$
\sum\limits_{t=0}^k (k+1 -t) h^1(\Cal O_C (t)) \leq (k+1) (k+2) g/2
\leq kd^2/4 + l.t.
\tag e
$$
Property (4) implies
$$
h^2 (\Cal O_X) = p_g (Y) \leq \frac{d^3}{24} + l.t.
\tag f
$$
So combining (a) $\sim$ (f) and applying Proposition 2, we have
$$
{\binom{k+5}{5}} = h^0 (\Cal O_{\P_5} (k)) \leq k^3 d/6 + d^2/2 +
d^3/6k + kd^2/4 + d^3/24 + l.t. \text{\qed}
$$
\enddemo

\proclaim{Proposition 5}  If $H^2K, HK^2 >0$, then $fHK^2 < .\ \ (d-
\delta) H^2 K$ holds for $X$ with $s(X) > f$, where $\delta$ is any
constant.
\endproclaim

\demo{Proof}  We may assume $f> 0$.  Assuming the contrary that $fHK^2
> (d-\delta)H^2 K$.  Multiplying by $H^2K$, {\rm applying GHIT (Property (2)),
and cancelling} $HK^2$, we have
$$
fH^2K > (d-\delta) d.
$$
The adjunction formula and genus bound give
$$
d^2/s + d(s-4)\geq 2g-2 > 2d+(d-\delta)d/f
$$
This is
$$
d< fs + f(f-6+\delta_1) + f^2 (f-6 +\delta_1)/(s-f),
\tag g 
$$
where $\delta_1 = \delta/f$.

By Proposition 4, (g) is false except for finitely many $s \leq s_0$.
For each $s\leq s_0$, there are only finitely many $d$ such that (g)
is true.  \qed
\enddemo

\proclaim{Proposition 6}  if $2-a-2b > 0$, then $a\Delta_1 + b\Delta_2 <.\ \
K^3$.
\endproclaim

\demo{Proof}  By Proposition 1 (ii), it suffices to show 
$$
\align
a[(d-15)H^2K - 6HK^2] + b[(6d-70) d+ (2d - 51)
H^2K - 12 HK^2]\\
<. \ \ 8d(d-15) + 2 (d-33)H^2K - 14HK^2.
\endalign
$$
By Proposition 5 and Proposition 1 (i), it suffices to show 
$$
(6b-8)HK^2 + (14-6a-12b)HK^2 <.\ \ (2-a-2b) d H^2 K.
$$
If $2-a-2b>0$, then this is 
$$
HK^2 (6-6a-6b)/(2-a-2b) <.\ \  d H^2 K .   
\tag h
$$
Let $f = (6 -6a-6b)/(2-a-2b)$.  Then Proposition 5 implies our claim for
those 3-folds $X$ with $s(X) > f$.

For $X$ with $s(X) \leq f$, the claim follows from Proposition 2. \qed
\enddemo

\subheading{Proof of the Theorem}

\demo{1.  Claim}  $y<.\ \ - x + 2 + \epsilon$, for any $\epsilon > 0$.
\enddemo

\demo{Proof}  This is $K^3 - \Delta_2 <.\ \ - K^3 + (2+\epsilon) ( K^3 - 
\Delta_1).$
\ \i.e.  $(1+\frac{2}{\epsilon}) \Delta_1 - \frac{1}{\epsilon} \Delta_2
<.\ \  K^3$.

We apply Proposition 6 with $a = 1+\frac{2}{\epsilon}$ and $b= - \frac{1}{\epsilon}$.
\enddemo

\demo{2.  Claim} $-\alpha x + 1 + \alpha - \epsilon <.\ \ y$, for any 
$\alpha > 1,\ \  \epsilon > 0$.
\enddemo

\demo{Proof}  This is 
$$
(1- (1+\alpha)/\epsilon) \Delta_1 + (1/\epsilon) \Delta_2 <.\ \  K^3
$$
Again, we take $a=1 - (1+\alpha)/\epsilon$ and $b = 1/\epsilon$ in 
Proposition 6.  Putting Claims 1 and 2 together, we have all limit points
lie on the line $x+y = 2$.

Similarly, we can show that $1-\epsilon <.\ \ x$ and $-\epsilon <.\ \ y$, for 
any $\epsilon > 0$.  Therefore there are no limit point with $x < 1$ or 
$y< 0$.  \qed
\enddemo

\remark{Remark}  In [C1] we showed that the limit points of $\{(c_1^3/c_1c_2, 
c_3/c_1c_2) (X)\; |\; X$ is the dependency locus of $\Cal O_{\P_5}^r \to
\opluss\limits_{i=1}^{r+1} \Cal O_{\P_5} (a_i)\}$ is the line segment
connecting the points $(1,1)$ and $(17/12, 7/12)$.
\endremark

\vfill\eject
\centerline{\bf References}
\ss
\roster
\item"{[BOSS]}"  R. Braun, G. Ottaviani, M. Schneider, and F. O. Schreyer,
Boundedness for 3-folds in $\P^5$.

\item"{[C1]}"  M. Chang, On the Chern numbers of surfaces and 3-folds of 
codimension 2, Tokyo J. of Math (to appear).

\item"{[C2]}" \lll, On zone SCI, (preprint)

\item"{[C3]}"  \lll, Appendix to ``On the geography of threefolds" by Liu,
Tohoku Math. J. (to appear).

\item"{[H]}"  B. Hunt, complex manifold geography in dimension 2 and 3, J. 
Differential Geometry, {\bf 30}, (1989), 51-153.

\item"{[L]}"  X. Liu, On the geography of threefolds, Tohoku Math. J. 
(to appear).
\endroster

\enddocument